\documentclass[12]{amsart}
\usepackage{mathrsfs}
\usepackage{amsfonts}
\usepackage{amssymb}
\newtheorem{thm}{Theorem}

\newtheorem{cor}[thm]{Corollary}
\newtheorem{lem}[thm]{Lemma}
\newtheorem{claim}[thm]{Claim}

\newtheorem{que}[thm]{Question}
\newtheorem{prop}[thm]{Proposition}
\newtheorem{exa}[thm]{Example}

\theoremstyle{definition}
\newtheorem{defn}[thm]{Definition}
\usepackage{amsmath,amscd}
\usepackage[all]{xy}
\newtheorem{rem}[thm]{Remark}

\begin{document}
 \title[Rings finitely generated left ideals are left annihilators of an element]{On rings whose finitely generated left ideals are left annihilators of an element}
\author[Yang]{Xiande Yang}
\date{\today}
\maketitle
 \centerline{\tiny { Department of Mathematics
and Statistics, University of New Brunswick, Fredericton, Canada E3B
5A3}}
 \centerline{\tiny Email address: xiande.yang@gmail.com}

\begin{abstract}
An associative ring $R$ with identity is left pseudo-morphic if for
every $a$$\in$$R$, there exists $b$$\in$$R$ such that $Ra=l_R(b)$.
If, in addition, $l_R(a)=Rb$, then $R$ is called left morphic. $R$
is  morphic if it is both left and right morphic. We characterize
left pseudo-morphic rings; identify the cases a (left) pseudo
morphic ring is (left) quasi-morphic, morphic, Quasi-Frobenius, von
Neumann regular, etc.;  correct two results in a book and a paper;
and completely determine when the trivial extension of a commutative
domain is morphic which positively answered a question in a paper.
\medskip



\end{abstract}

 \baselineskip=20pt
 \bigskip
\section{Introduction}
For an  associative ring $R$ with identity, Ehrlich \cite[Theorem
1(I)]{E76} proved that $R$ is unit regular iff $R$ is von Neumann
regular and for any $a\in R$, $\frac{R}{Ra}\cong l_R(a)$. Based on
this fact, Nicholson and S\'{a}nchez Campos \cite[Lemma 1]{NSC04}
defined a ring $R$ to be {\it left morphic} if for every element
$a\in R$ there exists $b\in R$ such that $\frac{R}{Ra}\cong l_R(a)$,
or equivalently,  $Ra=l_R(b)$ and $l_R(a)=Rb$. Left morphic rings
include unit regular rings, one-sided principal ideal artinian
rings, and some extensions of these rings \cite{D07, E76, LZ07,
LZ09, LZ10, NSC04, NSC042}. In a von Neumann regular ring, for any
$a\in R$, there exist $b, c\in R$ such that $Ra=l_R(b)$ and
$l_R(a)=Rc$. Camillo and Nicholson \cite{CN} called a ring $R$  {\it
left quasi-morphic} if for every element $a\in R$ there exist $b,
c\in R$ such that $Ra=l_R(b)$ and $l_R(a)=Rc$. Zhu and Ding
\cite{ZD} named a ring $R$ {\it left generalized morphic} if
$l_R(a)$ is principal for every $a\in R$. The right analogs are
defined similarly. $R$ is {\it morphic, quasi-morphic, or
generalized morphic} if $R$ is both left and right morphic,
quasi-morphic, or generalized morphic respectively. If we let
$\mathcal {LA}:=\{l_R(r)\text{ }| \text{ } r\in R\}$, $\mathcal
{LP}:=\{Rr\text{ }| \text{ } r\in R\}$ and $\mathcal
{LMP}=\{(a,b)\in R\times R\text{ }| (a,b) \text{ is a left morphic
pair, i.e., } l_R(a)=Rb \text{ and } l_R(b)=Ra\}$. Then we can say
$R$ is left generalized morphic  if  $\mathcal {LA}\subseteq
\mathcal {LP}$; $R$ is left quasi-morphic  if $\mathcal
{LA}=\mathcal {LP}$; and $R$ is left morphic if $\mathcal
{LA}=\mathcal {LP}$ and  the projection $p:\text{ }\mathcal
{LMP}\rightarrow R$  with $p(a,b)=a$ is surjective. There is a class
left--- rings with $\mathcal{LP}\subseteq \mathcal{LA}$. Yohe
\cite[Theorem 1]{Y68} called  $R$ a {\it left elemental annihilator
ring} ( l.e.a.r. for short) if every left ideal of $R$ is a left
annihilator of a single element of $R$ and proved that a commutative
ring is an l.e.a.r. iff it is a direct sum of artinian principal
local rings; Luedeman \cite[Theorem 5.2]{L70} and Yue Chi Ming
\cite[Theorem 2]{YCM70} proved that a nonsingular or semiprime ring
is an l.e.a.r. iff it is semisimple artinian; and Jaegermann and
Krempa \cite[Theorem 4.2]{JK72} characterized that left and right
elemental annihilator rings are exactly the direct sum  of matrix
rings over artinian principal left and right ideal rings. These
important rings satisfy the condition $\mathcal{LP}\subseteq
\mathcal{LA}$. The class of rings whose finitely generated one-sided
ideals are annihilator of an element include von Neumann regular
rings. Many papers are about when these rings are von Neumann
regular. For example, \cite[Theorem 1]{YCM02} is one. However, it
has a mistake. Motivated by these situations, we define left pseudo
morphic rings as the class of rings with $\mathcal{LP}\subseteq
\mathcal{LA}$ (Definition \ref{defn:1}).

In section 2, we characterize that left pseudo-morphic rings are
exactly rings whose finitely generated left ideals are left
annihilators of an element (Theorem \ref{thm:1}); give properties of
left pseudo-morphic rings (Theorem \ref{cor:03}); and disclose that
under some conditions a (left) pseudo-morphic ring is (left) quasi
morphic,  morphic, von Neumann regular, or Quasi Frobenius, and etc.
(Theorem  \ref{cor:04} and Theorem \ref{thm:5}).

In section 3, we correct minor mistakes  in \cite[Theorem
7.5(B)]{F99} (Theorem \ref{thm:2}) and  in \cite[Theorem 1]{YCM02}
(Theorem \ref{thm:3}) which essentially focused on the question when
a pseudo-morphic ring is von Neumann regular,

In section 4, we focus on commutative pseudo-morphic rings (they are
morphic at the same time). Specifically,  we focus on the trivial
extension $R\ltimes M$(see definition in section 4). Many papers are
on morphic properties of trivial extensions of a ring $R$ where $R$
is unit regular, strongly regular, one-sided perfect, or a domain.
Among them there are three for the case that $R$ is a domain.  Chen
and Zhou \cite[Theorem 14, Corollary 15]{CZ05} proved that ${\mathbb
Z}\ltimes \frac{{\mathbb Q}}{{\mathbb Z}}$ is strongly morphic  and
${\mathbb Z}\ltimes M$ is morphic iff $M\cong \frac{{\mathbb
Q}}{{\mathbb Z}}$ ($R$ is {\it strongly left morphic} if all matrix
rings over it are left morphic \cite{NSC05} and it is {\it strongly
morphic} if all matrix rings over it are morphic \cite{CZ05}). Lee
and Zhou \cite[Theorem 14]{LZ07} proved that when $R$ is a {\rm PID}
with classical quotient field $Q$, $R\ltimes \frac{Q}{R}$ is
strongly morphic and $R\ltimes M$ is morphic iff $M\cong
\frac{Q}{{R}}$. They further proved that if $R$ is a {\rm UFD}, then
$R\ltimes \frac{Q}{R}$ is morphic iff $R$ is a {\rm PID}. Recently,
Diesl, Dorsey,  and McGovern \cite[Theorem 4.13, Corollary 4.10,
Theorem 4.15, Question 1]{DDM09} characterized when $R\ltimes M$ is
left morphic for a domain $R$, showed that  $R\ltimes \frac{Q}{R}$
is morphic for $R$ a commutative B\'{e}zout domain, proved that if
$R$ is a commutative elementary divisor domain (see the definition
in Section 4) such that $R\ltimes M$ morphic, then $R\ltimes M$ is
strongly morphic (in particular, $R\ltimes \frac{Q}{R}$ is strongly
morphic), and asked whether or not for a commutative B\'ezout domain
$R$, $R\ltimes M$ is morphic iff $M\cong \frac{Q}{R}$. In Section 4,
we affirmatively answer this question (Theorem \ref{thm:6}) and
point out that when $R$ is an elementary divisor domain (need not be
commutative), $R\ltimes \frac{Q}{R}$ is strongly morphic
(Proposition \ref{thm:7}).

As usual, $J(R)$ and $U(R)$ denote Jacobson radical and the unit
group of $R$ respectively, and $l_R(a)$ and $r_R(a)$ (or $l(a)$ and
$r(a)$ if no confusion) are left and right annihilator ideals of $a$
in $R$.

\section{Left pseudo-morphic rings}

\begin{lem} Let $R$ be a ring and $a\in R$. Then the following are
equivalent.
\begin{enumerate}\item \label{item:1} There exists $b\in R$ such that $\frac{R}{Ra}\cong l(b)$.
\item \label{item:2} There exists $ c\in R$ such that $Ra=l(c)$, $Rc=l(b)$.
\item \label{item:3} There exists $ c\in R$ such that $Ra=l(c)$, $Rc\cong l(b)$.\end{enumerate}
\end{lem}
\begin{proof} ``$(\ref{item:1}) \Rightarrow$ (\ref{item:2})". Suppose   $\theta: \frac{R}{Ra}\rightarrow
l(b)$ is the isomorphism with $\theta(1+Ra)=c$. Then the following
diagram commutes where $\alpha$ and $\beta$ are canonical
epimorphisms and $\eta(1+l(c))=c$.
 \begin{center}  $\xymatrix{
                                                    R\ar[dr]^{\beta} \ar[r]^{\alpha}          & \frac{R}{Ra}\ar[r]^{\theta \text{  }}                &l(b)=Rc       \\
                                                                                              & \frac{R}{l(c)}\ar[ur]^{\eta}                &
  }$\end{center} Then we  get the $Ra=l(c)$ and $Rc=l(b)$.

  ``$(\ref{item:3}) \Rightarrow$ (\ref{item:1})".   $\frac{R}{Ra}=\frac{R}{l(c)}\cong Rc\cong l(b)$.
\end{proof}

\begin{defn}\label{defn:1} An element $a$ in a ring $R$ is called {\it left pseudo-morphic } if $a$ satisfies one of above equivalent conditions.
The ring $R$ is {\it left pseudo-morphic} if every  element of $R$
is left pseudo-morphic. Similarly we can define a  ring to be {\it
right pseudo-morphic}. $R$ is {\it pseudo-morphic} if it is both
left and right pseudo-morphic.
\end{defn}

Let $\mathcal {LA}:=\{l_R(r)\text{ }| \text{ } r\in R\}$, $\mathcal
{LP}:=\{Rr\text{ }| \text{ } r\in R\}$, and $\mathcal
{LMP}=\{(a,b)\in R\times R\text{ }| (a,b) \text{ is a left morphic
pair, i.e., } l_R(a)=Rb \text{ and } l_R(b)=Ra\}$. Then we can say
$R$ is left generalized morphic  if  $\mathcal {LA}\subseteq
\mathcal {LP}$; $R$ is left quasi-morphic  if $\mathcal
{LA}=\mathcal {LP}$; and $R$ is left morphic if $\mathcal
{LA}=\mathcal {LP}$ and  the projection $p:\text{ }\mathcal
{LMP}\rightarrow R$  with $p(a,b)=a$ is surjective.
\begin{cor} \label{cor:01}Let $R$ be a ring. Then $R$ is left pseudo-morphic if and
only if $Ra=l(b)$ for each $a\in R$ or, equivalently,
$\mathcal{LP}\subseteq \mathcal{LA}$.
\end{cor}
\begin{prop}\label{cor:02} Let $a$ be left pseudo-morphic in the ring $R$. Then the following hold.
\begin{enumerate}\item $ua$ and $au$ are left pseudo-morphic for any $u\in U(R)$.
\item If $r(a)=0$, then $Ra=R$. In particulr, $a$ is von Neumann regular and thus quasi-morphic.
\end{enumerate}
\end{prop}
\begin{proof}
$(1)$. By above lemma, we can suppose $Ra=l(b)$ and $Rb=l(c)$. Then
$Rau=l(b)u=l(u^{-1}b)$, $Rua=l(b)$, $Ru^{-1}b=l(c)$.

$(2)$. Let $Ra=l(b)$. Then $0=b\in r(a)=0$. So $Ra=R$. Thus, $a$ is
von Neumann regular with $a=ara=ea=af$ for some $r\in R$,  $e=ar$,
and $f=ra$. Thus, $Ra=l(1-f)$, $l(a)=R(1-e)$, $aR=r(1-e)$, and
$r(a)=(1-f)R$.
\end{proof}

\begin{thm}\label{thm:1} Let  $R$ be a ring. Then the following hold.
\begin{enumerate}\item $R$ is left generalized morphic iff any finite intersection $\cap_{i=1}^{n}l_R(a_i)=Ra$ for some $a\in
R$.
\item  $R$ is left pseudo-morphic iff any finite sum $\sum_{i=1}^{n}Ra_i=l_R(a)$ for some $a\in
R$.
\end{enumerate}\end{thm}
\begin{proof} $(1).$ Let $l_R(a_i)=Rb_i$ with $i=1,2$ and
$l_R(b_1a_2)=Rc$. Then $l_R(a_1)\cap\l_R(a_2)=Rcb_1$ since every
$x\in l_R(a_1)\cap\l_R(a_2)$ is of the form $x=r_1b_1=r_2b_2$ with
$r_1={r}^{'}c$ for some $r_{i}$ and $r^{'}$ in $R$, $i=1,2$.

$(2).$ Let $Ra_i=l_R(b_i)$ with $i=1,2$ and $Ra_2b_1=l_R(c)$. Then
$Ra_1+Ra_2=l_R(b_1c)$ since $xb_1c=0$ implies $xb_1=r_2a_2b_1$ for
some $r_2\in R$ and hence $x=r_2a_2+r_1a_1$ for some $r_1 \in R$.
\end{proof}

\begin{thm}\label{cor:03}If  $R$ is  a left pseudo-morphic ring, then the following hold.
\begin{enumerate}\item \label{item:4} $R$ satisfies  left annihilator conditions for finitely generated left ideals,
i.e., for any finitely generated left ideal $I$, $lr(I)=I$. \item
\label{item:5} $R$ is right principally injective, i.e., every
homomorphism from a principal right ideal to $R_R$ lifting to an
element of $R$.  In
 particular,
\begin{enumerate}\item $R$ is a right {\rm C2} ring, i.e., each right ideal isomorphic to a direct summand of $R_R$ is a direct summand. \item $J(R)=Z_r=\{a\in R \text{
}|\text{ } r(a) \text{ is essential in } R_R\}$.\item $R$ is right
mininjective, i.e. each homomorphism from a simple right ideal to
$R_R$ lifting to an element of $R$. So if $kR$ is simple, then $Rk$
is simple and  thus ${\rm S}_r\subseteq {\rm S}_l$ where $S_l$ is
the left socle and $S_r$ is the right socle of $R$.\item
\label{item:9} $R$ has {\rm ACC} on right annihilators iff $R$ is
right noetherian iff $R$ is right artianian.\end{enumerate}
\item \label{item:6} $R$ is directly finite iff $r(a)=0$ implies
$a\in U(R)$  for all $a\in R$ iff $r(a)=0$ implies $l(a)=0$ for all
$a\in R$.
\item \label{item:7} If $R_R$ has finite  Goldie dimension $Gdim(R)=m$, then $R$ is semilocal and has stable range one.
\item \label{item:8} If $R$ has $ACC$ on right annihilators $r(a)$
with $a\in R$, then $R$ is right perfect.

\end{enumerate}\end{thm}
\begin{proof}  $(\ref{item:4})$. By Theorem \ref{thm:1}, we can suppose that $I=\sum_{i=1}^{n}Ra_i=l_R(b)$. Then we have
$lr(I)=l_{R}r_{R}(\sum_{i=1}^{n}Ra_i)=l_Rr_Rl_R(b)=l_R(b)=\sum_{i=1}^{n}Ra_i=I$.

$(\ref{item:5})$. By $(\ref{item:4})$, for any $a\in R$,
$lr(a)=lr(Ra)=Ra$. So $R$ is right principally injective by
\cite[Theorem 1]{IN54}. Then we have $(a)$ by \cite[Proposition
5.10]{NY03}, $(b)$ by \cite[Theorem 5.14]{NY03}, $(c)$ by
\cite[Theorem 2.21]{NY03}, and $(d)$ by \cite[Proposition
5.15]{NY03}.

$(\ref{item:6})$.  Suppose $R$ is directly finite. By Proposition
\ref{cor:02}, for any $a\in R$, $r(a)=0$ implies $Ra=R$. So $ba=1$
for some $b\in R$. Thus, $a\in U(R)$ and $l(a)=0$.

 Suppose that $r(a)=0$ implies  $l(a)=0$ for all $a\in R$. If $ba=1$, then $r(a)=0$ and $aba=a$. So $l(a)=0$. Then
$(ab-1)a=0$ implies $ab=1$. Thus $R$ is directly finite.

Note this is essentially the same as \cite[Theorem 5.13]{NY03} if we
consider every element as a homomorphism.

$(\ref{item:7})$. Since $Gdim(R)=m$,  $R$ is directly finite. So, by
$(\ref{item:6})$, for any $a\in R$, $r(a)$=0 implies $a\in U(R)$.
Thus, by \cite[Theorem 5, statement on p.204]{CD93}, $R$ is
semilocal and has stable range one.

$(\ref{item:8})$. Let $Ra_1=l(b_1)\supseteq
Ra_2=l(b_2)\supseteq\ldots$ Then $r(a_1)\subseteq
r(a_2)\subseteq\ldots$ So there exists $n$ such that
$r(a_n)=r(a_{n+1})$, i.e. $r(Ra_n)=r(Ra_{n+1})$. So
$rl(b_n)=rl(b_{n+1})$, i.e. $Ra_n=Ra_{n+1}$. By \cite[Theorem
28.4(e)]{AF92}, $R$ is right perfect.
\end{proof}
\begin{cor} \label{cor:3}If $R$ is a left quasi-morphic ring, then the following hold.
\begin{enumerate} \item \cite[Thoerem 2, Corollary 4]{CNW08} Intersection of finitely many principal left
ideals of $R$ is a principal left ideal and  sum of finitely many
principal left ideals of $R$ is a principal left ideal. Thus, $R$ is
left B\'{e}zout and the set of principal left ideals $\mathcal {LP}$
forms a lattice with $Ra\vee Rb=Ra+Rb$ and $Ra\wedge Rb=Ra\cap Rb$.
\item \cite[Corollary 5]{CNW08} $R$ is left coherent.\item  If $R$ has {\rm ACC} on   principal left
ideals, then $R$ is a left noetherian principal left ideal
ring.\end{enumerate}
\end{cor}
\begin{proof}
$(1)$. By Theorem \ref{thm:1}(1),
$\cap_{i=1}^{n}Ra_i=\cap_{i=1}^{n}l_R({x_i}) =Rx$ for some $x,
x_i\in R$.  By Theorem \ref{thm:1}(2), $\sum_{i=1}^{n}Ra_i=l_R(x)
(=Ra)$ for some $a,x\in R$. The left is clear.

$(3)$. For any finitely generated left ideal $I$, $I=Ra$ since it is
left B\'{e}zout. Let $\scriptsize\xymatrix{0\ar[r]&{\rm
Ker}\theta\ar[r]^{i}&R^n\ar[r]^{\theta}&I=Ra\ar[r]&0 }\normalsize$
be any short exact sequence. Let ${\rm
Ker}\theta=\oplus_{j=1}^{n}I_j$, $\alpha_j=(0,\ldots, 0, 1,0\ldots
0)$ with the jth coordinate  1 , $K_j=R\alpha_j$ and
$\theta(\alpha_j)=r_ja$. Then
$\scriptsize\xymatrix{0\ar[r]&I_j\ar[r]^{i}&K_j\ar[r]^{{\theta|}_{K_j}}&Rr_ja\ar[r]&0
}\normalsize$ is exact. Since $R$ is left quasi-morphic,
$I_j=l(r_ja)=Rb$ for some $b\in R$. Hence ${\rm Ker}\theta$ is
finitely generated, i.e., $R$ is left coherent.

$(4)$. For any left ideal $I$ in $R$, the set of principal left
ideals in $I$ contains a maximal one $Ra$. Since $R$ is left
B\'{e}zout, $I=Ra$. Hence, $R$ is left noetherian.
\end{proof}
 The following example shows a left pseudo
morphic ring need not be right pseudo-morphic.

\begin{exa}(Nicholson and Yousif called this example of ring  Bj\"{o}rk example, see \cite{NY03}) Let $F$ be a field with a  homomorphism $\sigma: F\rightarrow F$ such that $\sigma(F)\not=F$ and let $R=\frac{F[x, \sigma]}{<x^2>}$ with $xa=\sigma(a)x$. Then $R$ is left morphic by \cite[Example 8]{NSC04}. In particular, it is left pseudo-morphic.
However, it is not right pseudo-morphic since $xR=\sigma(F)x$ is not
a right annihilator of a single element. (Suppose $xR=r_R(a+bx)$
with $a,b\in F$. By direct computation, $a=0$ and $b\not=0$. But
then $r_R(bx)=Fx\not=\sigma(F)x=xR$ which is a contradiction.)
\end{exa}

For a local ring $R$ with $J(R)$ nilpotent, $R$ is left generalized
morphic iff $R$ is left quasi-morphic iff $R$ is left morphic iff
$R$ is an artinian principal left ideal ring(see \cite[Lemma 9]{CN},
\cite[Proposition 2.8]{ZD}, and \cite[Theorem 9]{NSC04}). This
explains why the artinian local rings ${\mathbb Z}_4$ and ${\mathbb
Z}_4\ltimes 2{\mathbb Z}_4$ are morphic  but $R={\mathbb Z}_4\ltimes
{\mathbb Z}_4\cong\frac{{\mathbb Z}_4[x]}{<x^2>}$ and
$T=\frac{F[x_1, x_2,\ldots]}{ <x_ix_j, i,j=1,2,\ldots>}$ with $F$ a
field are neither left generalized morphic  nor left pseudo-morphic
(In fact, $l_R(2x)=2R+xR$ is not principal, $2xR=l_R(2)\cap l_R(x)$
is not an annihilator of a single element, and $T$ is  not
artinian). The following example shows the condition ``local" is
necessary and points out that even under the condition ``artinian" a
left generalized morphic ring need not be a left pseudo-morphic.

\begin{exa}Let
$R=\scriptsize\left(\begin{array}{ll}\mathbb{Z}_2&0\\\mathbb{Z}_2&\mathbb{Z}_2\end{array}\right)\normalsize.$
Then $R$ has 5 idempotents and 2 units which are left morphic. The
only undetermined element  is
$\alpha=\scriptsize\left(\begin{array}{ll}0&0\\1&0\end{array}\right)\normalsize$.
By direct computation,
$l_R(\alpha)=R\scriptsize\left(\begin{array}{ll}1&0\\1&0\end{array}\right)\normalsize$
which is principal but $R\alpha=\left\{0,
\scriptsize\left(\begin{array}{ll}0&0\\1&0\end{array}\right)\normalsize\right\}$
$=l_R\left\{\scriptsize\left(\begin{array}{ll}1&0\\1&0\end{array}\right)\normalsize
\right\}\cap
l_R\scriptsize\left\{\left(\begin{array}{ll}0&0\\0&1\end{array}\right)\normalsize
\right\}$ is not a left annihilator of one element. So $R$ is left
generalized morphic but  not left pseudo-morphic. In fact, $R$ is
not left morphic (hence, not left quasi-morphic) can also be
determined by \cite[Proposition 18 ]{NSC04}.
\end{exa}

Every domain which is not a division ring is generalized morphic but
not pseudo-morphic (hence, not quasi-morphic). It is well known that
every unit regular ring is morphic and  every von Neumann regular
ring which is not unit regular is quasi-morphic but not morphic. So
we have the relations: \noindent $\{\text{Left generalized-morphic
Rings}\}\supset \{\text{Left Quasi-morphic Rings}\} \supset
\{\text{Left Morphic Rings}\},$ \noindent $\{\text{Left pseudo
morphic rings}\}\supseteq \{\text{Left Quasi-morphic Rings}\}
\supset \{\text{Left Morphic Rings}\}$.

\begin{que}Is a left pseudo-morphic ring left quasi-morphic? Or
equivalently, is a left pseudo-morphic ring left generalized
morphic?
\end{que}
A ring $R$ is a {\it left elemental annihilator ring} ( l.e.a.r. for
short) if every left ideal of $R$ is a left annihilator of a single
element of $R$ (see \cite{Y68}). A ring $R$ is {\it left
Ikeda-Nakayama} if $r(I_1\cap I_2)=r(I_1)+r(I_2)$ for any left
ideals $I_1$ and $I_2$ in $R$ (see \cite{CNY00}). Note,  generally,
$r(I_1\cap I_2)\supseteq r(I_1)+r(I_2)$ for any subset $I_1, I_2\in
R$.  The following results shows the relations of (left)
pseudo-morphic rings with  (left) quasi-morphic rings and other
important rings.
\begin{thm}\label{cor:04} Let $R$ be a ring.
\begin{enumerate}
\item The following are equivalent.
\begin{enumerate}
\item $R$ is left quasi-morphic.
\item $R$ is left pseudo-morphic, left coherent, and left B\'{e}zout.
\end{enumerate}
\item The following are equivalent.
 \begin{enumerate}
 \item $R$ is pseudo-morphic.
 \item $R$ is quasi-morphic.
 \item \cite[Corollary 4]{CN} $R$ is morphic if, in addition,  $R$ is commutative.

 \noindent In above cases, $R$ is an {\rm IF} ring ( $R$ is an {\rm IF} ring if all injective left or right $R$-modules are flat),
  $r\left(\cap_{j=1}^{n} I_j\right)=\sum_{j=1}^{n}r(I_j)$ for  finitely generated left ideals $I_j$, and $l\left(\cap_{j=1}^{n}
I_j\right)$ $=\sum_{j=1}^{n}l(I^{'}_j)$ for  finitely generated
right ideals $I^{'}_j$.
 \end{enumerate}

 \item The following are equivalent for a pseudo-morphic ring $R$.

 \begin{enumerate}\item $R$ is left noetherian.
\item $R$ is an l.e.a.r and an r.e.a.r.
\item $R=\oplus {\mathbb M}_{n_{i}}(R_i)$ with $R_i$ artinian local principal left and right ideal ring.
\item $R$ is a Quasi-Frobenius ring.
\item For any left or right module $M$ of $R$, $M$ is injective if and only
if $M$ is projective iff $M$ is flat. The class of these modules are
closed under direct sum and direct product.
\item $R$ is a right noetherian ring.
\item $R$ is a left artinian ring.
\item $R$ is a right artinian ring.

In above cases,  $R$ is a dual ring, i.e., $lr(I)=I$ and $rl(T)=T$
for all left ideal $I$ and right ideal $T$; the matrix ring
${\mathbb M}_n(R)$ is a strongly clean ring and an Ikeda-Nakayama
ring for any integer $n>0$; and $R$ is semisimple artinian if, in
addition, $R$ is semiprime (or nonsingular).
\end{enumerate}
\end{enumerate}\end{thm}
\begin{proof} $(1)$. ``$(a)\Rightarrow (b)$. By Corollary \ref{cor:3}.

 ``$(b)\Rightarrow (a)$. Given any
$x\in R$, by the exact sequence $0\rightarrow l_R(x)\rightarrow
R\rightarrow Rx\rightarrow0$, we can assume
$l_R(x)=\sum_{i=1}^nRa_i$ since $R$ is left coherent. Then
$l_R(x)=\sum_{i=1}^nRa_i=Ra$ since $R$ is left B\'{e}zout.

$(2)$. ``$(a)\Rightarrow (b)$".  If $R$ is pseudo-morphic, we show
that for any $a\in R$, $l(a)$ and $r(a)$ are principal.  Since $R$
is pseudo-morphic, we can assume that $aR=r(b)$ and $Rb=l(c)$. Then
$l(a)=l(aR)=lr(b)=lr(Rb)=lrl(c)=l(c)=Rb$. Symmetrically, we can
prove that $r(a)$ is principal. So $R$ is quasi-morphic.

 ``$(b)\Leftrightarrow (c)$". \cite[Corollary 4]{CN}, i.e., commutative
 quasi-morphic rings are  morphic.

   If $R$ is quasi-morphic, then $R$ is
coherent by Corollary \ref{cor:3} and every finitely generated
one-sided ideal is an annihilator of an element by Theorem
\ref{thm:1}(2). Then $R$ is an {\rm IF} ring by \cite[Theorem
6.9]{F99}. By Theorem \ref{cor:03}(2), $R$ is left and right
principally injective. By Corollary \ref{cor:3}(3), $R$ is
B\'{e}zout.  Then by \cite[Theorem 1(ii)]{IN54},
$r\left(\cap_{j=1}^{n} I_j\right)=\sum_{j=1}^{n}r(I_j)$ and
$l\left(\cap_{j=1}^{n} I_j\right)=\sum_{j=1}^{n}l(I^{'}_j)$ for any
finitely generated left ideals $I_j$ and finitely generated right
ideals $I^{'}_j$.

$(3)$.  ``$(a)\Rightarrow (b)$".  By Theorem \ref{cor:03}(2.d) and
the fact that $R$ is Bezout, we get that $R$ is a left artinian
principal left ideal ring. So every left ideal is a left annihilator
of an element. Since $R$ is right quasi-morphic, $aR=r(a^{'})$ for
all $a\in R$. By \cite[Theorem 4.2(iii)]{JK72}, $R$ is an l.e.a.r.
and an r.e.a.r.

``$(b)\Rightarrow (c)$". By \cite[Theorem 4.2(vii)]{JK72}.

``$(c)\Rightarrow (d)$". Since $R$ is an artinian and principally
injective. So $R$ is self-injective.

``$(d\Rightarrow (e)$". By $(2)$, $R$ is an {\rm IF} ring, i.e.,
every injective module is flat,  thus projective. But $R$ is self
injective and artinian, so every projective module is injective.

``$(e)\Rightarrow (a)$". By \cite[Proposition 18.13]{AF92}.

The equivalence of $(f)$, $(g)$, and $(h)$ to others can be proved
similarly.

Now, in above cases, $R$ is self-injective and thus ${\mathbb
M}_n(R)$ is Ikeda-Nakayama by \cite[Theorem 7]{CNY00}. Every
artinian ring is strongly ${\pi}$-regular, hence, strongly clean
\cite[Proposition 1]{N99}. $R$ is dual by $(3.b)$. If, in addition,
$R$ is a semiprime or nonsingluar, then $R$ is semisimple artinian
by \cite[Theorem 5.2]{L70} and \cite[Theorem 2]{YCM70}.
\end{proof}

\begin{rem}We proved above results  by the  old result
\cite[Theorem 4.2]{JK72}. In fact, \cite[Theorem 19, Lemma 18]{CN}
are equivalent to Theorem \ref{cor:04}.3.c  by  Corollary
\ref{cor:3}(3).
\end{rem}

\begin{thm}\label{thm:5}Let $R$ be a reduced ring (i.e., $R$ contains no nonzero nilpotent element)and $n$ any positive integer. Then the following are equivalent.
\begin{enumerate}
\item $R$ is left pseudo-morphic
\item $R$ is left quasi-morphic.
\item $R$ is left morphic.
\item $\frac{R[x]}{<x^n>}$ is left pseudo-morphic.
\item $\frac{R[x]}{<x^n>}$ is left quasi-morphic.
\item $\frac{R[x]}{<x^n>}$ is left morphic.
\item $R$ is regular.
\item $R$ is unit regular.
\item $R$ is strongly regular.
\end{enumerate}
change "left" into "right" or delete ``left" in the above items, the
results are still equivalent to $(1)$.
\end{thm}
\begin{proof}$(1)\Rightarrow (9).$ For any $a\in R$, let $Ra=l_R(b)=r_R(b)$
with $b\in R$. Note  $x\in l_R(a+b)=r_R(a+b)\Rightarrow ax=-bx
\Rightarrow axa=0$ $\Rightarrow (bx)^2=0$  $\Rightarrow bx=0=xb$
$\Rightarrow x=ra$ for some $r\in R$. Now $0=-bx=ax=ara \Rightarrow
x^2=rara=0\Rightarrow x=0$. So $l_R(a+b)=r_R(a+b)=0$. Let
$R(a+b)=l_R(c)$. Then $(a+b)c=0$. So $c\in r_R(a+b)=l_R(a+b)=0$. We
get $R(a+b)=l_R(0)=R$. So $1=r(a+b)$ for some $r\in R$. Hence
$a=r(a+b)a=ra^2$. So $R$ is strongly regular.

  $(9)\Leftrightarrow (8)\Leftrightarrow (7)$ since idempotents are central in $R$:
$(ex-exe)^2=(xe-exe)^2=0$$ \Rightarrow ex=xe$ for any $x \text{ and
} e^2=e\in R$.

$(8)\Rightarrow (6)$ by \cite[Theorem 11]{LZ09}.

$(6)\Rightarrow(5)\Rightarrow(4)$. It is clear from the definitions.

$(4)\Rightarrow (1).$ Here, we prove a more generalized case that
for any ring $R$, if $\frac{R[x]}{<x^n>}$ is left pseudo-morphic,
then $R$ is left pseudo-morphic: Let $T=\frac{R[x]}{<x^n>}$. For any
$a\in R$, $ax^{n-1}\in T$ and $Tax^{n-1}=Rax^{n-1}$. So there exists
$b=r_0+r_1x+\ldots+r_{n-1}x^{n-1}$ with $r_{i}\in R$, $i=0,\ldots,
(n-1)$, such that $Tax^{n-1}=Rax^{n-1}=l_T(b)$. So $Ra\subseteq
l_R(r_0)$. For any $r\in l_R(r_0)$, $rx^{n-1}\in l_T(b)$, i.e.,
$r\in Ra$. So $Ra= l_R(r_0)$, i.e., $R$ is left pseudo-morphic.

$(8)\Rightarrow (3)$: This is  \cite[Theroem 1(I)]{E76}.

$(3)\Rightarrow (2)\Rightarrow (1)$: It is clear.
\end{proof}

A ring $R$ is {\it reversible} if any $ab=0$ implies $ba=0$ and $R$
is {\it symmetric} if any $abc=0$ implies $acb=bac=0$. For a ring,
the property of being reduced or symmetric always implies that of
reversible. In above theorem, the condition ``reduced" can not be
replaced by ``reversible" or ``symmetry"  since  $R=\frac{{\mathbb
Z}_{2}[x]}{<x^2>}$ is symmetric, left pseudo-morphic by above
results but $R$ is not even von Neumann regular.

\section{corrections on rings Pseudo-morphic  implying  von Neumann regular}
In this section, we correct two minor mistakes in a book and a
paper. As Theorem \ref{thm:1} showed, rings with every finitely
generated left ideal left annihilator of an element are exactly left
pseudo-morphic rings. To determine if a semiprime pseudo-morphic
ring is von Neumann regular, Yue Chi Ming \cite{YCM02} got the
following result.
 \begin{claim}\label{claim:2}\cite[Theorem
1]{YCM02} The following conditions are
equivalent:\begin{enumerate}\item $A$ is von Neumann regular.
\item A is a semiprime  ring whose finitely generated one-sided ideals are annihilators of an element of $A$ (this is equivalent to that $R$ is semiprime pseudo-morphic). \item $A$ is a semiprime ring such that every
finitely generated left ideal is the left annihilator of an element
of $A$ and every principal right ideal of $A$ is the right
annihilator of an element of $A$ (this is still equivalent to that
$R$ is semiprime pseudo-morphic).\end{enumerate}\end{claim}

In the proof of Claim \ref{claim:2}, the author used the following
result.

\begin{claim}\label{claim:1} \cite[Theorem 7.5(B)]{F99}  A semiprime right
B\'{e}zout ring $R$ is right semi-hereditary and a finite product of
prime right B\'{e}zout rings each isomorphic to a full $n\times n$
matrix ring ${\mathbb M}_n(F)$ over a right Ore domain $F$ for
various integers $n\geq 1$.
\end{claim}

However, there exists  a semiprimitive (thus, semiprime) right
B\'{e}zout left quasi-morphic (thus left pseudo-morphic) ring which
is neither semi-hereditary nor von Neumann regular.
\begin{exa} \label{exa:1} For a subring $C$ of a ring $D$, let $\mathcal
{R}[D,C]=\{(d_1, \ldots, d_n, c, c, \ldots) : d_i\in D, c\in C,
n\geq1\}$ with componentwise addition and multiplication. Then
$\mathcal {R}[D,C]$ is a ring. Let $S=\mathcal {R}[{\mathbb
M}_n(R),\frac{R[x]}{<x^n>}]$ with $R$ von Neumann regular. As
\cite[Theorem 9]{LZ09} showed $S$ is semiprimitive (thus
semiprime)and quasi-morphic (hence right B\'{e}zout \cite{CNW08}).
By  Claim \ref{claim:1}, $S$ is semi-hereditary. Thus every finitely
generated left ideal of $S$ is projective and so generated by an
idempotent. Hence, $S$ is von Neumann regular. However, as
\cite[Theorem 9]{LZ09} showed that $S$ is not von Neumann regular
(in fact, the homomorphic image of $(x, x, \ldots)$ in
$\frac{R[x]}{<x^n>}$ is $x$ which is not von Neumann regular).
\end{exa}
According to \cite[Corollary 4.2 and Corollary 4.4]{R67} and the
fact that right Ore domain is exactly right Goldie domain,  Claim
\ref{claim:1}  can be changed to the following.

\begin{thm}\cite[Corollary 4.2 and Corollary 4.4]{R67} \label{thm:2} A ring $R$ is semiprime, right B\'{e}zout, and right Goldie iff
$R=\sum_{i=1}^{k}{\mathbb M}_{n_{i}}(R_i)$ where ${\mathbb
M}_{n_{i}}(R_i)$ is prime right B\'{e}zout with $R_i$ a right Ore
domain. In particular, $R$ is semi-hereditary.
\end{thm}

The proof of Claim \ref{claim:2} depends on Claim \ref{claim:1}. So
if we add the condition that $R$ is right Goldie, then we can change
Claim \ref{claim:2} to the following.
\begin{thm} \label{thm:3} Let $R$ be a right Goldie ring. Then $R$ is a semiprime  pseudo-morphic (or
equivalently, quasi-morphic by Theorem \ref{cor:04}) iff $R$ is von
Neumann regular iff  $R$ is semisimple artinian.\end{thm}
\begin{proof}If $A$ is pseudo-morphic, then by Theorem \ref{cor:04}, $A$ is B\'{e}zout.
Then by  Theorem \ref{thm:2}, $A$ is left semihereditary. Hence
every principal left ideal is projective and generated by an
idempotent. Thus $R$ is von Neumann regular.

If $R$ is von Neumann regular, then $R= \sum_{i=1}^{k}{\mathbb
M}_{n_{i}}(R_i)$ where ${\mathbb M}_{n_{i}}(R_i)$ is prime right
B\'{e}zout with $R_i$ a right Ore domain by Theorem \ref{thm:2}.
Then $R$ is a semisiple artinian since matrix ring over a domain is
von Neumann regular iff the domain is a division ring.

When $R$ is semisimple artinian, it is clear that $R$ is seimprime
pseudo-morphic.
\end{proof}
We require $A$ to be right Goldie because the author's proof used
Theorem \ref{thm:2}. In fact, we can change \cite[Theorem 1]{YCM02}
to the following.
\begin{thm}Let $R$ be a ring. Then $R$ is von Neumann regular if and
only if $R$ is left semihereditary (or left p.p.) and  right pseudo
morphic (or left principally injective) iff $R$ is right
semihereditary (or right p.p.) and  left pseudo-morphic (or right
principally injective).
\end{thm}
\begin{proof} We only prove that a left p.p. right pseudo-morphic ring
is von Neumann regular. $R$ is left principally injective by
Corollary \ref{cor:03} and  $Ra$ is projective for any $a\in R$.
Then $Ra$ is a direct summand of $R$ by \cite[Corollary 5.11]{NY03}.
So $a$ is von Neumann regular.
\end{proof}

\section{Morphic trivial extension of a commutative domain}
As Theorem  \ref{cor:04}(2) shown,  a commutative  pseudo-morphic
ring is morphic. In this section, we focus on commutative pseudo
morphic rings.   First, we give necessary conditions for some left
generalized morphic extensions and left pseudo-morphic extensions;
second, we completely determine when the trivial extension of a
commutative domain is morphic, this affirmatively answered a
question in \cite{DDM09}.

The {\it trivial extension} of a ring $R$ by an $R$-bimodule $M$ is
$R\ltimes M=\{(r,m): r\in R \text{ and } m\in M\}$ with
componentwise addition and multiplication $(r_{1}, m_{1})(r_{2},
m_{2})=(r_{1}r_{2}, r_{1}m_{2}+m_{1}r_{2})$. Note that $R\ltimes
M\cong\{\bigl(\begin{smallmatrix}r&m\\0&r\end{smallmatrix}\bigr):
r\in R, m\in M\}$ and if $M=R$ then $R\ltimes R\cong
\frac{R[x]}{<x^2>}$.
\begin{prop}\label{prop:12} Let $R$ and $S$ be rings and $_RV_S$ be an $R$-$S$ bimodule and $C=\scriptsize\left(
\begin{array}{cc}
R&V\\
0&S
\end{array}
\right)\normalsize$. Then the following hold.
\begin{enumerate}\item If $C$ is left generalized morphic, then both $R$ and
$S$ are left generalized morphic.\item If $C$ is left pseudo
morphic, then $S$ is left pseudo-morphic. If $C$ is right left
pseudo-morphic, then $R$ is right pseudo-morphic.\end{enumerate}
\end{prop}
\begin{proof} $(1)$. Let $C$ be left generalized morphic. Let $\alpha=\scriptsize\left(
\begin{array}{cc}
a&0\\
0&0
\end{array}
\right)\normalsize \in C$. Then $\scriptsize\left(
\begin{array}{cc}
r&v\\
0&s
\end{array}
\right)\normalsize \alpha=\scriptsize\left(
\begin{array}{cc}
ra&0\\
0&0
\end{array}
\right)$. So $l_C(\alpha)= \scriptsize\left(
\begin{array}{cc}
l_R(a)&V\\
0&S
\end{array}
\right)\normalsize$. Since $l_C(\alpha)$ is principal, there exists
an element $\beta=\scriptsize\left(
\begin{array}{cc}
x&y\\
0&z
\end{array}
\right)\normalsize\in C$ such that $\scriptsize\left(
\begin{array}{cc}
R&V\\
0&S
\end{array}
\right)\normalsize \beta=\scriptsize\left(
\begin{array}{cc}
Rx&Ry+Vz\\
0&Sz
\end{array}
\right)\normalsize = l_C(\alpha).$ Hence, $l_R(a)=Rx$ is principal,
i.e., $a$ is left generalized morphic. So $R$ is left generalized
morphic ring.  Let $\gamma=\scriptsize\left(
\begin{array}{cc}
0&0\\
0&a
\end{array}
\right)\normalsize$. By a similar argument we get that $S$ is left
generalized morphic ring.

$(2)$. If $C$ is left pseudo-morphic, let
$\alpha=\scriptsize\left(\begin{array}{ll}0&0\\0&a\end{array}\right)\normalsize$.
If $C$ is right pseudo-morphic, let
$\alpha=\scriptsize\left(\begin{array}{ll}a&0\\0&0\end{array}\right)\normalsize$.
The by the similar argument of above, we can prove the results.
\end{proof}

\begin{cor}\label{cor:13} Suppose $e^2=e\in R$ such that $(1-e)Re=0$. Then the following hold.
\begin{enumerate}
 \item If $R$ is left generalized morphic, then so are $eRe$ and $(1-e)R(1-e)$.\item If $R$ is left pseudo-morphic, so is $(1-e)R(1-e)$. If
$R$ is right pseudo-morphic, so is $eRe$.\end{enumerate}
\end{cor}
\begin{proof}
By the Pierece decomposition $R \cong \scriptsize\left(
\begin{array}{cc}
eRe&eR(1-e)\\
(1-e)Re&(1-e)R(1-e)
\end{array}
\right)\normalsize$.
\end{proof}

\begin{prop}Let $T=R\ltimes M$ be the trivial extension of $R$ by the bimodule $_RM_R$. If $T$ is  left
generalized morphic, then  $R$ is left generalized
morphic.\end{prop}
\begin{proof}Suppose $T$ is left generalized morphic. For any $a \in R$, let  $\alpha=(a,0)\in T$. Then
$l_T(\alpha)=l_R(a)\ltimes l_{M}(a)=T\beta=\{(rx,rn+mx): r\in R,
m\in M\}$ with $\beta=(x,n)\in T$. So $l_R(a)=Rx$ is left
generalized morphic.
\end{proof}

So when we do morphic trivial extensions, we require the base ring
be generalized morphic. So we focus on trivial extension over a
commutative domain $R$.
 \begin{lem}\label{lem:21}Let $R$ be a commutative domain with classical quotient field $Q$ and let $R\overline{\frac{1}{a}}$
 and $R\overline{\frac{1}{b}}$ be submodules of $\frac{Q}{R}$ with $a,b\in R$. Then the following hold.
 \begin{enumerate}\item  $R\overline{\frac{1}{b}}\leq R\overline{\frac{1}{a}}$ iff $b | a$.
 \item $R\overline{\frac{1}{a}}\cong R\overline{\frac{1}{b}}$ iff $R\overline{\frac{1}{a}}= R\overline{\frac{1}{b}}$ with
 $a=rb$, $r\in
 U(R)$.
 \item If $R$ is B\'{e}zout, then every cyclic submodule of $\frac{Q}{R}$ is of the form
$R\overline{\frac{1}{c}}$, $c\in R$.\end{enumerate}\end{lem}
 \begin{proof} $(1)$. Suppose $R\overline{\frac{1}{b}}\leq
 R\overline{\frac{1}{a}}$. Then $\overline{\frac{1}{b}}=r\overline{\frac{1}{a}}$. Hence
 $a-br=kab$ for some $k\in R$. So $b(ka+r)=a$.

 $(2)$. Note that $R\overline{\frac{1}{a}}\cong R\overline{\frac{1}{b}}$ implies
 $Rb=l_R(\overline{\frac{1}{b}})=l_R(\overline{\frac{1}{a}})=Ra$.
 Hence $a=rb$ for some $r\in U(R)$.

 $(3)$. Since $R$ is B\'ezout, for any $r, s\in R$ with $r\not=0$ or $s\not=0$, there exists some $d\in R$ such that $Rr+Rs=Rd$. Hence $r=r^{'}d$ and $s=s^{'}d$
for some $r^{'}, s^{'}\in R$. So $Rr^{'}+Rs^{'}=R$, i.e., $(r^{'},
s^{'})$ is a unimodular pair in $R$. Now we can assume every cyclic
submodule of $\frac{Q}{R}$ is of the form $R\overline{\frac{r}{c}}$
for some unimodular pair $(r, c)$ in $R$ with $c\not=0$. Hence
$rr^{'}+cc^{'}=1$ for some $r^{'}, c^{'}\in R$. So in $Q$ we have
$\frac{rr^{'}}{c}+\frac{cc^{'}}{c}=\frac{1}{c}$. Therefore, in
$\frac{Q}{R}$ we have $\overline{\frac{1}{c}}=\overline
{\frac{rr^{'}}{c}}$, i.e.,
$R\overline{\frac{r}{c}}=R\overline{\frac{1}{c}}$.
 \end{proof}
 \begin{prop}\label{prop:1}Let $R$ be a commutative B\'ezout domain with classical quotient field
$Q$. Then  the following hold.

\begin{enumerate}\item The set of  finitely generated submodules of $\frac{Q}{R}$ is
$\mathscr{S}(\frac{Q}{R})=\{R\overline{\frac{1}{a}}:a\in R\}$. \item
$\mathscr{S}(\frac{Q}{R})$ forms a lattice. \item There are no
isomorphic members in $\mathscr{S}(\frac{Q}{R})$ except the member
itself.\end{enumerate}
\end{prop}
\begin{proof}$(1)$. By Lemma \ref{lem:21}(3), the set of finitely generated submodules is
$\mathscr{S}(\frac{Q}{R})=\{R\overline{\frac{1}{a}}:a\in R\}$.

$(2)$. Given any $a, b\in R$. Suppose $Ra+Rb=Rd$. Define
$R\overline{\frac{1}{a}}\wedge
R\overline{\frac{1}{b}}=R\overline{\frac{1}{d}}$ and
$R\overline{\frac{1}{a}}\vee
R\overline{\frac{1}{b}}=R\overline{\frac{d}{ab}}$. Then
$\mathscr{S}(\frac{Q}{R})$ forms a lattice.

$(3)$. By Lemma \ref{lem:21}(2), there are no isomorphic members in
$\mathscr{S}$ except they are equal.
\end{proof}
For our convenience, we include two results of \cite{DDM09} here.
\begin{lem}\label{lem:4}(\cite[Corollary 4.10]{DDM09}) If $R$ is
a commutative B\'ezout domain with classical ring of quotient $Q$,
then $R\ltimes \frac{Q}{R}$ is morphic.\end{lem}

\begin{lem}\label{lem:5}( \cite[Corollary 4.14]{DDM09}) Suppose that  $R$ is
a  domain which is not a division ring and that $M$ is an
$R$-bimodule such that $R\ltimes M$ is morphic. Then the map
$\mathcal {F}(mR)=l_R(m)$ is an inclusion-reversing bijection from
the set $\mathscr{S}(M)=\{mR: m\in M\}$ to the set $\{Ra: 0\not=a\in
R\}$.\end{lem}
\begin{thm} \label{thm:6} Let $R$ be a commutative domain with classical quotient field $Q$ and let $M$ be an $R$-module.
Then the following hold. \begin{enumerate}\item If $R$ is a field,
then $R\ltimes M$ is morphic iff $M=0$ or $M\cong R$.
\item If $R$ is not a field, then  $R\ltimes M$ is morphic if and only
if $R$ is B\'ezout and $M\cong \frac{Q}{R}$.\end{enumerate}
\end{thm}
\begin{proof} Note the trivial extension $R\ltimes M$ is commutative if and only
if $R$ is commutative. So being left morphic, right morphic, or
morphic are the same for the ring $R\ltimes M$. So in the following
proof we do not mention left or right properties such as
 B\'ezout and so on. And here we use $Rn < Rm$ to mean $Rn$ is a
 proper submodule of  $Rm$.

$(1)$. This is a part of \cite[Proposition 11]{LZ07}(see also
\cite[Corollary 4.12]{DDM09}).

$(2)$.   $``\Leftarrow"$. This is Lemma \ref{lem:4}.

 $``\Rightarrow"$. Suppose that
$R\ltimes M$ is morphic, then  $R$ and $M$ are B\'ezout by
\cite[Proposition 2.4]{DDM09}. For any $Rn < Rm\leq M$, we can let
$l_{R}(n)=Rb$ and $l_{R}(m)=Ra$ by Lemma \ref{lem:5}. Then $Rb>Ra$
by Lemma \ref{lem:5}. Now we have that $Rm\cong \frac{R}{Ra}\cong
R\overline{\frac{1}{a}}$ and $Rn\cong \frac{R}{Rb}\cong
R\overline{\frac{1}{b}}< R\overline{\frac{1}{a}}$. By Lemma
\ref{lem:21}.1, $a=bb^{'}$ with some $b^{'}\in R$. Notice that the
isomorphism $Rm\cong R\overline{\frac{1}{a}}$ can be defined by
$f_{(m,a)}: m\mapsto \overline{\frac{1}{a}}$ (the proof is
standard). Clearly $f^{-1}_{(m,a)}(\overline{\frac{1}{b}})=b^{'}m$.
Suppose that $n=r_{n}m$. Then
$l_{R}(n)=l_{R}(r_{n}m)=Rb=l_{R}(\overline{\frac{1}{b}})=l_{R}(b^{'}m)$.
Then $Rn\cong Rb^{'}m$. By Lemma \ref{lem:5},
$\mathscr{S}(M)\cong\mathscr{S}(\frac{Q}{R})$. By Proposition
\ref{prop:1}, $\frac{Q}{R}$ has no isomorphic finitely generated
submodules except identical ones.  So $M$ has no isomorphic finitely
generated submodules except identical ones. Hence $Rn= Rb^{'}m$. So
we have a commutative diagram where $i$ is the inclusion and
$f_{(b^{'}m, b)}$ is defined similar to  $f_{(m, a)}$.
\begin{center}$\xymatrix{
  Rm  \ar[r]^{f_{(m,a)}} & R\overline{\frac{1}{a}}  \\
  Rn=Rb^{'}m \ar[u]_{i}\ar[r]^{ \text{  }f_{(b^{'}m, b)}} & R\overline{\frac{1}{b}}.\ar[u]^{i}
  }$\end{center}

By the above commutative diagram, we can suppose that in the
following diagram every triangle diagram and parallelogram  diagram
commute except the bottom one.

    \begin{center}  $\xymatrix{
                                                   & Rm \ar[rr]^{f_{(m,a)}}          &                 & R\overline{\frac{1}{a}}       \\
                                                   & Rb^{'}m\ar[u]_{i_{3}}\ar[rr]^{f_{(b^{'}m,\text{ } b)}}&                 &R\overline{ \text{  }\frac{1}{b}\text{  }}\ar[u]_{i_{4}}  \\
  Rc^{'}m\ar[uur]^{i_{1}}\ar[ur]_{i_{5}}\ar[rr]^{f_{(c^{'}m, \text{ },c)}}                    &                     &R\overline{\frac{1}{c}}.\ar[uur]_{i_{2}}\ar[ur]_{i_{6}}&
  }$\end{center}

 Note that $i_{4}f_{(b^{'}m,
b)}i_{5}= f_{(m, a)}i_{3}i_{5}=f_{(m, a)}i_{1}=i_{2}f_{(c^{'}m,
c)}=i_{4}i_{6}f_{(c^{'}m, c)}$. Since $i_{4}$ is the inclusion, we
have $f_{(b^{'}m, b)}i_{5}=i_{6}f_{(c^{'}m, b)}$. So the bottom
parallelogram diagram commutes.

Now we have two directed system $\{Rm: m\in M \}$ and
$\{R\overline{\frac{1}{a}}: a\in R\}$ with inclusions as module
homomorphisms. Furthermore, $M=\underrightarrow{\text{lim}}\{Rm:
m\in M \}$ and
$\frac{Q}{R}=\underrightarrow{\text{lim}}\{R\overline{\frac{1}{a}}:
a\in R\}$ since both $M$ and $\frac{Q}{R}$ are B\'ezout. By above
commutative diagrams
 we have the following commutative
diagram.

\begin{center}$\xymatrix{M=\underrightarrow{\text{lim}}\{Rm: m\in M \}\ar@/_3mm/[rrr]^{g}&&&\frac{Q}{R}=\underrightarrow{\text{lim}}\{R\overline{\frac{1}{a}}:a\in R\}\ar@/_5mm/[lll]_{h}\\
&\vdots&\vdots&\\
&Rm\ar[uul]_{\beta_{m}}\ar[u]\ar[r]^{f_{(m,a)}}&R\overline{\frac{1}{a}}\ar[u]\ar[uur]^{\alpha_{a}}&\\
&Rb^{'}m\ar[uuul]^{\beta_{b^{'}m}}\ar[u]\ar[r]^{f_{(b^{'}m, b)}}&R\overline{\frac{1}{b}}\ar[u]\ar[uuur]_{\alpha_{b}}&\\
&\vdots\ar[u]&\vdots\ar[u]&
 }$\end{center}

Since $f_{(m,a)}$ are isomorphisms, by the universal property of
direct limits, we get that $g$ and $h$ are isomorphisms. So $M\cong
\frac{Q}{R}$.
\end{proof}
By above theorem, we  generalize \cite[Theorem 14]{LZ07} as
following.
\begin{cor}\label{cor:1}  Let $R$ be a  {\rm UFD} or a commutative noetherian domain  with classical quotient field $Q\not=R$ and $M$ be an $R$-bimodule. Then
$R\ltimes M$ is morphic iff $R$ is a {\rm PID} and $M\cong
\frac{Q}{R}$.
\end{cor}
\begin{proof} $``\Rightarrow".$ By Theorem \ref{thm:6}, $M\cong
\frac{Q}{R}$ and  $R$ is B\'ezout. It is well-known that when $R$ is
a B\'ezout domain, then $R$ is a ${\rm PID}$ iff $R$ is a {\rm UFD}
iff $R$ is noetherian.

$``\Leftarrow".$ By Theorem \ref{thm:6} or Lemma \ref{lem:4}.
\end{proof}

 \begin{rem}We point out  there are many B\'ezout domains which are not {\rm
PID} (see \cite{C68}).\end{rem}

An associative ring $R$ with unit is an elementary divisor ring if
every matrix over $R$ has a diagonal reduction, i.e., for every
matrix $A$ over $R$, there exist invertible matrices $P$ and $Q$
over $R$  such that $PAQ={\rm diag}(d_1, d_2, \cdots, d_{r},  0,
\cdots, 0)$ is a diagonal matrix and $Rd_{i+1}R\subseteq d_{i}R\cap
Rd_{i}$ for $i=1, \cdots, r-1$ (see \cite{GZ99}). In \cite[Theorem
4.15]{DDM09}, the authors proved that if $R$ is a commutative
elementary divisor domain and $M$ is a bimodule (in fact, $M\cong
\frac{Q}{R}$ by Theorem \ref{thm:6}) such that $R\ltimes M$ is
morphic, then $R\ltimes M$ is strongly morphic. In fact, it can be
generalized  to a non-commutative domain.
\begin{prop} \label{thm:7} If $R$ is an elementary divisor domain (need not be commutative) and $M$ is a
bimodule such that $R\ltimes M$ is morphic, then $R\ltimes M$ is
strongly morphic.
\end{prop}
\begin{proof}The proof is similar to that  in \cite{DDM09} except we need to change the bimodule $_RM_R$ into $_{R-R^{opp}}M$ where $R^{opp}$ is the opposite ring of $R$.\end{proof}

\section*{Acknowledgements}
The research was done under the supervision of Professor Colin
Ingalls. Thanks to Professor Carl Faith and Professor J. Chris
Robson for their kind help during preparing the manuscript.
Appreciation also goes to Professor Hugh Thomas for the help and
comments on an earlier version. The author is partially supported by
ACEnet of Canada as a postdoctoral fellow.

\end{document}